\documentclass{amsart}
\usepackage{amssymb}
\usepackage{amsmath}
\usepackage{amsthm}

\newtheorem{theorem}{Theorem}[section]
\newtheorem{proposition}{Proposition}[section]
\newtheorem{lemma}[theorem]{Lemma}
\theoremstyle{definition}
\newtheorem{definition}[theorem]{Definition}

\theoremstyle{remark}
\newtheorem{remark}[theorem]{Remark}
\numberwithin{equation}{section}

\def\H{\mathcal H}
\def\Soc{\operatorname{Soc}}
\def\Top{\operatorname{Top}}
\def\Rad{\operatorname{Rad}}
\def\Hom{\operatorname{Hom}}
\def\rank{\operatorname{rank}}

\begin{document}
\date{}
\title{Proof of the modular branching rule for cyclotomic Hecke algebras}
\author{Susumu Ariki}
\address{Research Institute for Mathematical Sciences,\\
         Kyoto University, Kyoto 606-8502. Japan}
\email{ariki@kurims.kyoto-u.ac.jp}

\begin{abstract}
We prove the modular branching rule of the cyclotomic Hecke algebras, 
which has remained open. 
\end{abstract}
\maketitle

\section{Introduction}

Let $F$ be an algebraically closed field. 
The cyclotomic Hecke algebra $\H_n=\H_n(\underline v,q)$ of type 
$G(m,1,n)$ is the $F$-algebra introduced in \cite{AK} and \cite{BM}. 
This is a cellular algebra in the sense of 
Graham and Lehrer, and the cell module theory of this 
algebra is nothing but the Specht module theory developed 
by Dipper, James and Mathas \cite{DJM1}. 
\footnote{Specht module theory 
for Hecke algebras was initiated by Dipper and James, 
and the Specht module theory we use here is its 
generalization to the cyclotomic Hecke algebras.} 
The Specht modules are parametrized by $m$-tuples of partitions 
$\lambda=(\lambda^{(1)},\dots,\lambda^{(m)})$ and denoted by 
$S^\lambda$. Each $S^\lambda$ has an invariant symmetric 
bilinear form, and we denote by $D^\lambda$ the module 
obtained from $S^\lambda$ by factoring out the radical of 
the invariant form. Then nonzero $D^\lambda$'s form a complete 
set of irreducible $\H_n$-modules. 

If we set $m=1$, $\H_n$ is the Hecke algebra of 
type $A$. If we further set $q=1$, then $\H_n$ is the group 
algebra of the symmetric group $S_n$. Kleshchev studied 
$\Soc(\operatorname{Res}^{S_n}_{S_{n-1}}(D^\lambda))$ in a series of 
papers \cite{Kl1} to \cite{Kl4}, 
and obtained an explicit rule for describing the socle. 
This is called the modular branching rule of the symmetric group. 
The method is to use modified lowering operators, and Brundan generalized 
this result to the Hecke algebra of type $A$ by the same method \cite{B}. 

Around the same time, motivated by conjectures and results by Lascoux, 
Leclerc and Thibon, a link between quantum groups of type 
$A^{(1)}_{e-1}$ and the Hecke algebra of type $A$ was found. 
In particular, they observed that the crystal rule of 
Misra and Miwa coincides with Kleshchev's rule for 
the modular branching \cite{LLT}. 

On the other hand, 
in solving the LLT conjecture on the decomposition numbers, 
I generalized the LLT conjecture to the graded 
dual of Grothendieck groups of the module categories of 
$\H_n$ with common parameters. With this interpretation, 
the action of Chevalley generators 
is given by refined restriction and induction functors, 
which are the 
$i$-restriction and $i$-induction functors. 
\footnote{The use of central elements in 
the symmetric group goes back to Robinson, which I learned from 
Leclerc, but the 
refined induction and restriction operators in this context 
were introduced by the author.} 
Further, by using Lusztig's canonical 
basis in the proof, it was natural for us to observe the existence 
of a crystal structure on the set 
$$
B=\bigsqcup_{n\ge0}\{\text{isoclasses of simple $\H_n$-modules}\}.
$$
\footnote{This was already mentioned in the form of its relationship with 
Kashiwara's lower crystal basis in \cite[p.807]{A1}.}
In this theory, 
which we call Fock space theory, 
we may identify the crystal with $\mathcal{KP}$ 
of those multipartitions for which 
$D^\lambda\ne0$. \footnote{We named these multipartitions Kleshchev 
multipartitions in \cite{AM}.} Its rigidity, namely independence of 
the characteristic of $F$, was first proved in \cite{AM}. 
The crystal is isomorphic to the $\mathfrak g(A^{(1)}_{e-1})$-crystal 
of an integrable highest weight module $L_v(\Lambda)$, where 
$e$ is the multiplicative order of the parameter $q\ne1$ and 
$\Lambda$ is determined by the parameters $\underline v$. 
For the overview of the Fock space theory, see \cite{Abook}. 

As many people in our field noticed, these works give a natural 
conjecture generalizing the
results of Kleshchev and Brundan on modular branching rules for the symmetric
groups and the Hecke algebras of type A; that is, we have a natural conjecture
for a modular branching rule for the cyclotomic Hecke algebras. 
Explicitly, this asserts that $\Soc(e_iD^\lambda)=D^{\tilde e_i\lambda}$, 
where $e_i$ is the $i$-restriction and 
$\tilde e_i$ is the Kashiwara operator of the crystal $\mathcal{KP}$. 

There was a progress toward this conjecture in Vazirani's thesis, 
which was later published as \cite{GV}. In the thesis, various facts which 
are necessary to show that $B$ has a crystal structure are proven, 
and they are used in \cite{G} to show that our $B$ is equipped 
with another crystal structure. This crystal structure is again 
isomorphic to the crystal of the same integrable highest module 
\cite[Theorem 14.3]{G}. 
In fact, the proof is carried out 
within the framework of my Fock space theory. 

On this occasion, I correct two of his announcements which are relevant to the
modular branching rule, as service to the mathematical community and to avoid
confusions.  In \cite{GV}, it is said: \lq\lq What we do not do in this paper
is to explicitly describe which irreducible representations occur in the socle
of the restriction.  This is done in \cite{G}, generalizing [Kv](=Kleshchev's
work) which describes the combinatorics of the branching rule for the symmetric
group explicitly in terms of $p$-regular partitions.\rq\rq However, in \cite{G}
one only finds such a result in terms of an abstractly defined crystal graph,
and no attempt is made to give an explicit description of the latter in terms
of partitions. Moreover, Grojnowski left completely untouched the problem of
matching up the standard labeling of simple modules coming from Specht module
theory with his labeling coming from the abstract crystal graph. So, contrary
to the announcement recorded in the note added in proof of Mathas' book
\cite[p.135]{Mbook}, no proof of the modular branching rule (even in the 
case of Brundan and Kleshchev's original modular
branching rule) is present in \cite{GV} or in \cite{G} 
(except in the case where
$q$ is not a root of unity which was treated in \cite{V}). 

As modular branching is used as the definition of their crystal, it is more
appropriate to see their theory as a method to label simple modules using a
crystal, rather than as a modular branching rule. \footnote{In fact, viewing
the theory this way, Brundan and Kleshchev were able to label simple modules of
the Hecke-Clifford algebra by using $\mathfrak g(A^{(2)}_{2l})$-crystal.} The
adjoint operation to modular branching is to take head of induction, and they
use this as the method to label simple $\H_n$-modules. This means that we need
to repeat the operation of taking the head of an induced module $n$ times to
compute a simple $\H_n$-module this way.

Let us examine in more detail how to compute the label of a given module,
and modular branching, by this method. 
Suppose that we are given a 
simple $\H_n$-module $V$ and that we have computed its 
character, namely its restriction to the commutative subalgebra 
generated by Jucy-Murphy elements. Then we can compute 
the character of $e_iV$. To know $\Soc(e_iV)$, 
we have to rewrite the character into summation over 
characters of simple $\H_{n-1}$-modules $V'$ and compute the values 
$\epsilon_i(V')$. Thus we are required to know the 
irreducible characters. 
\footnote{In modular representation theory, 
knowing irreducible characters is a hard problem. One may list 
the modular representation theory of the symmetric group, 
the Kazhdan-Lusztig and Lusztig conjectures, as examples. 
Note that knowing irreducible characters is equivalent to 
computing decomposition numbers.} 
The only way to compute the irreducible 
characters in the method is to construct the modules 
by taking head of induction as above. 
\footnote{Here, Specht module theory provides us with 
easier way to construct the simple modules, but it is still 
unrealistic to compute the irreducible characters by constructing 
simple modules.}
One can compute the character of an induced module, but 
we meet the same problem for computing its head. 
Thus, to compute the labeling or modular branching, 
the only way is to compute socle of restriction (or 
head of induction) explicitly. 

Finding the label $\lambda$ of a given module in 
Specht module theory is also not automatic, but we have more realistic 
chance for finding the label. For example, $\lambda$ is the 
minimal Kleshchev multipartition that satisfies 
$\Hom_{\H_n}(S^\lambda,V)\ne0$. Further, 
the original modular branching rule allows us 
to compute the socle of the restriction without 
computing the socle. It is also worth mentioning that our 
approach of using the Specht module theory is still the only 
alternative even for proving Brundan's result in type $A$. 
Thus, the importance of the Specht module theory 
could not be overestimated. 

The purpose of this paper is to prove the  
modular branching rule of cyclotomic Hecke algebras, 
which until now has remained 
open. It turns out that it 
is a direct consequence of the theorem on 
the canonical basis in the Fock space. 
\cite{GV} and \cite{G} contain new results also 
and we use two of them in the proof. 
\footnote{When writing this paper, I learned that 
Brundan had a very similar idea for the proof. 
He considered a similar problem in a different setting 
\cite[Theorem 4.4]{BK2}, and observed that 
the same strategy works in the present situation. 
I thank Brundan for the communication.}

\section{Preliminaries}

\begin{definition}
Let $R$ be a commutative ring, and let $v_1,\dots,v_m, q\in R$ 
be invertible elements. The cyclotomic Hecke algebra 
$\H_n(\underline v,q)$ is the $R$-algebra defined by 
the generators $T_0,\dots,T_{n-1}$ and the relations
\begin{gather*}
(T_0-v_1)\cdots(T_0-v_m)=0,\quad (T_i-q)(T_i+1)=0,\;\text{for}\;i\ge1,\\[5pt]
(T_0T_1)^2=(T_1T_0)^2,\\[5pt]
T_iT_{i+1}T_i=T_{i+1}T_iT_{i+1}, \;\text{for}\;i\ge1,\\[5pt]
T_iT_j=T_jT_i, \;\text{for}\;j\ge i+2.
\end{gather*}
\end{definition}

We write $\H_n$ for short. It is known that $\H_n$ is free of rank $m^nn!$ 
as an $R$-module. We define elements $L_1,\dots,L_n$ by 
$$
L_1=T_0,\quad L_{k+1}=q^{-1}T_kL_kT_k,\;\text{for}\;1\le k<n.
$$
They pairwise commute and the symmetric polynomials in $L_1,\dots,L_n$ 
are central elements of $\H_n$. 

The Specht module theory for $\H_n$ is developed by 
Dipper, James and Mathas \cite{DJM1}. Recall that 
the set of multipartitions, namely the set of $m$-tuples of partitions, of 
size $n$ is a poset whose partial order is the dominance 
order $\trianglerighteq$. 
Let $\lambda=(\lambda^{(1)},\dots,\lambda^{(m)})$ be a 
multipartition of size $n$. 
Then we can associate an $\H_n$-module $S^\lambda$ with $\lambda$, 
called a Specht module. $S^\lambda$ is free as an $R$-module. 
Further, each Specht module is equipped with 
an invariant symmetric bilinear form \cite[(3.28)]{DJM1}. 
Let ${\rm rad}\,S^\lambda$ be the radical of the invariant 
symmetric bilinear form, and 
we set $D^\lambda=S^\lambda/{\rm rad}\,S^\lambda$. 
We denote the projective cover of $D^\lambda$ by $P^\lambda$ 
when $D^\lambda\ne0$. 

\begin{theorem}[{\cite[Theorem 3.30]{DJM1}}]
Suppose that $R$ is a field. Then, 

\begin{itemize}
\item[(1)]
Nonzero $D^\lambda$ form 
a complete set of non-isomorphic simple $\H_n$-modules. Further, 
these modules are absolutely irreducible. 
\item[(2)]
Let $\lambda$ and $\mu$ be multipartitions of size $n$ and 
suppose that $D^\mu\ne 0$ and that $[S^\lambda:D^\mu]\ne 0$. 
Then $\lambda\trianglerighteq\mu$. Further, 
$[S^\lambda:D^\lambda]=1$. 
\end{itemize}
\end{theorem}

The projective cover $P^\mu$ has a Specht filtration 
$$
P^\mu=F_0\supset F_1\supset\cdots
$$
such that $F_0/F_1\simeq S^\mu$. This follows from the cellularity of $\H_n$. 

By the Morita-equivalence theorem of Dipper and Mathas \cite{DM}, we may 
assume that $v_i$ are powers of $q$ without loss of generality. 
In the rest of paper, we assume that $q$ is a primitive $e^{th}$ 
root of unity where $e\ge2$, and 
$v_i=q^{\gamma_i}$, for $\gamma_i\in\mathbb Z/e\mathbb Z$. 

\section{The Kashiwara crystal}

Let $A=(a_{ij})_{i,j\in I}$ be a generalized Cartan matrix, 
$\mathfrak g=\mathfrak g(A)$ the Kac-Moody Lie 
algebra associated with $A$. 
Let $(P,\Delta,P^{\rm v},\Delta^{\rm v})$ be the simply-connected 
root datum of $\mathfrak g$. 
We write $\alpha_i$ for simple roots, and $h_i$ for simple coroots. 
Thus, 
$P^{\rm v}$ is generated by $\{h_i\}_{i\in I}$ and $|I|-\rank(A)$ 
elements $\{d_s\}$ as a $\mathbb Z$-module. 

\begin{definition}
A $\mathfrak g$-crystal $B$ is a set endowed with 
\begin{itemize}
\item
$wt:B \longrightarrow P$,
\item
$\epsilon_i,\varphi_i: B \longrightarrow \mathbb Z\sqcup\{-\infty\}$,
\item
$\tilde e_i,\tilde f_i: B \longrightarrow B\sqcup\{0\}$,
\end{itemize}
such that the following properties are satisfied. 
\begin{itemize}
\item[(1)]
$\varphi_i(b)=\epsilon_i(b)+\langle h_i,wt(b)\rangle$.
\item[(2)]
If $b\in B$ is such that $\tilde e_ib\ne0$ then
$$
wt(\tilde e_ib)=wt(b)+\alpha_i,\;\;
\epsilon_i(\tilde e_ib)=\epsilon_i(b)-1,\;\;
\varphi_i(\tilde e_ib)=\varphi_i(b)+1.
$$
\item[(3)]
If $b\in B$ is such that $\tilde f_ib\ne0$ then
$$
wt(\tilde f_ib)=wt(b)-\alpha_i,\;\;
\epsilon_i(\tilde f_ib)=\epsilon_i(b)+1, \;\;
\varphi_i(\tilde f_ib)=\varphi_i(b)-1.
$$
\item[(4)]
For $b,b'\in B$, we have $b'=\tilde e_ib \Longleftrightarrow \tilde f_ib'=b$. 
\item[(5)]
If $b\in B$ is such that $\varphi_i(b)=-\infty$ then 
$\tilde e_i(b)=0$ and $\tilde f_i(b)=0$. 
\end{itemize}
\end{definition}

Let $U_v(\mathfrak g)$ be the quantized enveloping algebra and 
$L_v(\Lambda)$ an integrable highest weight 
$U_v(\mathfrak g)$-module. Then the lower crystal base $B(\Lambda)$ 
of $L_v(\Lambda)$ is a $\mathfrak g$-crystal. Further, the crystal 
$B(\Lambda)$ is semiregular. That is, 
$$
\epsilon_i(b)={\rm max}\{k\in\mathbb Z_{\ge0}|\tilde e_i^kb\ne0\},\quad
\varphi_i(b)={\rm max}\{k\in\mathbb Z_{\ge0}|\tilde f_i^kb\ne0\}.
$$

The module $L_v(\Lambda)$ has a distinguished basis, which is called 
Kashiwara's (lower) global basis or the Lusztig canonical basis. 
The basis elements are labelled by $B(\Lambda)$, and we denote them by 
$\{G_v(b)\}_{b\in B(\Lambda)}$. See \cite{HK} for example. 

The following lemma is taken from \cite[Lemma 12.1]{K2}. For the proof, follow 
the argument in \cite[Proposition 5.3.1]{K1} which is for the upper global basis. 

\begin{lemma}\label{Kashiwara}
Let $B(\Lambda)$ be the crystal of the integrable highest weight module 
$L_v(\Lambda)$. Then the following hold. 
\begin{itemize}
\item[(1)]
There exist Laurent polynomials $e_{bb'}^i(v)$ such that 
$$
e_iG_v(b)=[\varphi_i(b)+1]G_v(\tilde e_ib)+\sum_{b'}e_{bb'}^i(v)G_v(b'),
$$
where the sum is over $b'\in B(\Lambda)$ with 
$\varphi_j(b')\ge\varphi_j(b)+\langle h_j,\alpha_i\rangle$, for all $j$. 
\item[(2)]
There exist Laurent polynomials $f_{bb'}^i(v)$ such that 
$$
f_iG_v(b)=[\epsilon_i(b)+1]G_v(\tilde f_ib)+\sum_{b'}f_{bb'}^i(v)G_v(b'),
$$
where the sum is over $b'\in B(\Lambda)$ with
$\epsilon_j(b')\ge\epsilon_j(b)+\langle h_j,\alpha_i\rangle$, for all $j$. 
\end{itemize}
\end{lemma}

In this paper, we only use the affine Kac-Moody Lie algebra of type 
$A^{(1)}_{e-1}$, where $e$ is defined by the parameter $q$ as in the 
previous section. The crystal we use is 
the $A^{(1)}_{e-1}$-crystal $B(\Lambda)$, where 
$\Lambda=\sum_{i=1}^m \Lambda_{\gamma_i}$ and $\gamma_i$ are 
$v_i=q^{\gamma_i}$ as before. 

\section{Fock space theory}

The Fock space theory is explained in detail in \cite{Abook}. 
Let $\mathfrak g$ be the affine Kac-Moody Lie algebra of type $A^{(1)}_{e-1}$. 
In \cite{A1}, I introduced the combinatorial Fock space $\mathcal F(\Lambda)$. 
It is a based $\mathbb Q$-vector space whose basis is 
the set of all multipartitions $\mathcal P$. 
The weight $\Lambda$ defines a rule to 
color nodes of multipartitions with $e$ colors $\mathbb Z/e\mathbb Z$, 
and the coloring rule defines an integrable $\mathfrak g$-module 
structure on $\mathcal F(\Lambda)$. 
Its deformation $\mathcal F_v(\Lambda)$ 
becomes an integrable $U_v(\mathfrak g)$-module via the Hayashi action, 
and the crystal obtained from $\mathcal F_v(\Lambda)$ is $\mathcal P$. 
Let $W_i(\lambda)$ be the number of $i$-nodes in $\lambda$. 
Then by the definition of the Hayashi action, we have
\begin{equation*}
\begin{split}
wt(\lambda)(h_i)&
=\Lambda(h_i)+W_{i-1}(\lambda)-2W_i(\lambda)+W_{i+1}(\lambda), 
\;\text{for}\; 0\le i\le e-1,\\
wt(\lambda)(d)&
=\Lambda(d)-W_0(\lambda).
\end{split}
\end{equation*}

Recalling 
$\alpha_j(h_i)=a_{ij}$ and $\alpha_j(d)=\delta_{j0}$, this is equivalent to 
$$
wt(\lambda)=\Lambda-\sum_{j=0}^{e-1}W_j(\lambda)\alpha_j.
$$
Kashiwara operators $\tilde e_i$ and $\tilde f_i$ are defined by 
removing or adding a good $i$-node. As $\mathcal P$ is semiregular, 
$\epsilon_i$ and $\varphi_i$ are determined by $\tilde e_i$ and $\tilde f_i$. 
Then $(\mathcal P,\tilde e_i, \tilde f_i, wt, \epsilon_i, \varphi_i)$ 
is the crystal structure given on $\mathcal P$. 

The connected component of $\mathcal P$ that contains the empty multipartition 
$\emptyset$ is denoted by $\mathcal{KP}$, and we call 
multipartitions in $\mathcal{KP}$ Kleshchev multipartitions. 
The global basis $\{G_v(\lambda)\}_{\lambda\in\mathcal{KP}}$ is the basis 
of the $U_v(\mathfrak g)$-submodule generated by $\emptyset$, which is 
isomorphic to the irreducible highest weight $U_v(\mathfrak g)$-module 
$L_v(\Lambda)$. Similarly, 
the basis $\{G_v(\lambda)\}_{\lambda\in\mathcal{KP}}$ evaluated at $v=1$ 
is the basis of the $\mathfrak g$-submodule generated by $\emptyset$, 
which is isomorphic to the irreducible highest weight 
$\mathfrak g$-module $L(\Lambda)$. We denote 
$\{G_v(\lambda)\}_{\lambda\in\mathcal{KP}}$ evaluated at $v=1$ by 
$\{G(\lambda)\}_{\lambda\in\mathcal{KP}}$. 

Suppose that the ground ring $R$ of $\H_n$ 
is an algebraically closed field $F$ of characteristic $\ell$, 
and recall that $q$ is a primitive $e^{th}$ root of 
unity and $v_i=q^{\gamma_i}$, for $\gamma_i\in \mathbb Z/e\mathbb Z$. 
Let $\H_n\text{-}proj$ be the category of (finite dimensional) 
projective $\H_n$-modules. 
In \cite{A1}, 
I defined the $i$-restriction and the $i$-induction functors. 
Let us recall the definitions following \cite[13.6]{Abook}. 
Let $M$ be an $\H_n$-module. 
As symmetric polynomials in $L_1,\dots,L_n$ are central elements in $\H_n$, 
the simultaneous generalized eigenspace with respect to 
the symmetric polynomials in $L_1,\dots,L_n$ is again an $\H_n$-module. 
Let $c=\{c_1,\dots,c_n\}$ where $c_i\in q^{\mathbb Z}$, and denote 
by $P_c(M)$ the simultaneous generalized eigenspace 
which consists of $m\in M$ such that 
$$
(f(L_1,\dots,L_n)-f(c_1,\dots,c_n))^Nm=0, 
$$
for $N>\!>0$ and for symmetric polynomials $f$. Define 
$$
e_iM=\sum_c P_{c\setminus\{q^i\}}
(\operatorname{Res}^{\H_n}_{\H_{n-1}}(P_c(M))), \;\;
f_iM=\sum_c P_{c\cup\{q^i\}}
(\operatorname{Ind}^{\H_n}_{\H_{n-1}}(P_c(M))).
$$
$e_i$ is the $i$-restriction and $f_i$ is the $i$-induction. 
They are exact functors. 

Suppose that $R$ is a discrete valuation ring, $K$ its fraction field, and 
$M$ an $\H_n$-module which is torsionless as a $R$-module. Then 
we have $M\subset M\otimes_R K$, where $M\otimes_R K$ is an 
$\H_n\otimes_R K$-module, and 
the definitions of $e_i$ and $f_i$ make sense for $M$. 
Further, $e_i(M\otimes_R K)=(e_iM)\otimes_R K$ and 
$f_i(M\otimes_R K)=(f_iM)\otimes_R K$ hold. 

The following are main results of \cite{A1}. 
See sections (4.5), (4.6), 
Theorem 4.4 and Proposition 4.5 in \cite{A1}, or 
Theorem 12.5 and Proposition 13.41 in \cite{Abook}. 

\begin{theorem}\label{Fock space theory}
Let $K_0(\H_n\text{-}proj)$ be the Grothendieck group of $\H_n\text{-}proj$. 
Then 
\begin{itemize}
\item[(1)]
The action of $e_i$ and $f_i$ on 
$K(\Lambda)=\oplus_{n\ge0}K_0(\H_n\text{-}proj)$ satisfy the Serre 
relations, and extends to a $\mathfrak g$-module structure on 
$K(\Lambda)$. 
\item[(2)]
$K(\Lambda)$ is isomorphic to 
the integrable $\mathfrak g$-module $L(\Lambda)$. 
\item[(3)]
We have a unique injective $\mathfrak g$-module homomorphism 
$K(\Lambda)\rightarrow \mathcal F(\Lambda)$ which sends 
the highest weight vector $[P^\emptyset]$ 
to the empty multipartition $\emptyset$. 
\item[(4)]
Assume that the characteristic of $F$ is zero, and that 
$D^\lambda\ne0$. Then $[P^\lambda]$ 
maps to a basis element $G(\lambda')$, for 
some $\lambda'\in\mathcal{KP}$, and we have 
$$
G(\lambda')=\lambda+\;\text{(higher terms)}\;=
\sum_{\mu\trianglerighteq\lambda}d_{\mu\lambda}\mu,
$$
where $d_{\mu\lambda}=[S^\mu:D^\lambda]$, the decomposition numbers. 
\end{itemize}
\end{theorem}

Note that the existence of a crystal structure on the set 
$$
B=\bigsqcup_{n\ge0}\{\text{isoclasses of simple $\H_n$-modules}\}, 
$$
is clear from this theorem. 
That $\lambda'=\lambda$ is proved in \cite{A2}. 
In particular, $D^\lambda\ne0$ if and only if 
$\lambda\in\mathcal{KP}$ and we can identify $B$ with $\mathcal{KP}$. 

For each simple module $D^\lambda$, we have that any 
symmetric polynomial $f$ in $L_1,\dots,L_n$ acts as a 
scalar. Because of our assumption that $v_i$ are powers of $q$, 
the eigenvalues of $L_k$, for $1\le k\le n$, are powers of $q$. 
This is because they are powers of $q$ for Specht modules. 
Thus, we have a uniquely determined set 
$\{q^{i_1},\dots,q^{i_n}\}$ such that 
every symmetric polynomial $f(L_1,\dots,L_n)$ acts on $D^\lambda$ as 
the scalar 
$f(q^{i_1},\dots,q^{i_n})$. 
Observe that the symmetric polynomials act as scalars on $S^\lambda$ 
already, and we can describe the set $\{q^{i_1},\dots,q^{i_n}\}$ 
explicitly as follows. 
$$
|\{k\in[1,n]| q^{i_k}=q^i\}|=W_i(\lambda).
$$
This module theoretic interpretation of $W_i(\lambda)$ was 
used in \cite{A1}, and will be used in the next section. 

\section{Another crystal structure}

Grojnowski and Vazirani introduced another 
semiregular crystal structure on the set $B$. 
The $i$-restriction they use \cite[3.1]{GV} is precisely the one 
which I introduced in \cite{A1}. 
$f_i$ is left and right adjoint to $e_i$. As one can see from 
the definition of $f_i$ given before, the definition 
is in terms of generalized eigenspace of $L_n$. 
Grojnowski introduced another description of $f_i$ \cite[p.17]{G}. 
If one observes that the $i$-restriction gives Jordan block of $L_n$, 
this description of $f_i$ is quite natural and not surprising at all. 
However, the point is that Vazirani and Grojnowski systematically 
developed properties of my functors and this approach is 
more suitable to study the 
modular branching rule. 
The crystal structure may be defined as follows. 
$$
\tilde e_iD^\lambda=\Soc(e_iD^\lambda),\quad
\tilde f_iD^\lambda=\Top(f_iD^\lambda), \quad
wt(D^\lambda)=wt(\lambda).
$$

As the crystal we define is semiregular, $\epsilon_i$ and $\varphi_i$ are 
determined by $\tilde e_i$ and $\tilde f_i$. As is stated in the introduction, 
the following is proved in \cite[Theorem 12.3]{G}. 

\begin{theorem}\label{Brundan-Kleshchev}
Let $(B, \tilde e_i, \tilde f_i, wt, \epsilon_i, \varphi_i)$ 
be as above. Then $B$ is isomorphic to $B(\Lambda)$. 
\end{theorem}

Another result of Grojnowski and Vazirani \cite[Lemma 3.5]{GV} implies that 
we can detect $\tilde e_iD^\lambda$ on the Grothendieck group level. 

\begin{proposition}\label{Grothendieck level}
If $\tilde e_iD^\lambda\ne0$, 
$\tilde e_iD^\lambda\ne0$ is a unique composition factor 
$D^\mu$ of $e_iD^\lambda$ 
with $\epsilon_i(D^\mu)=\epsilon_i(D^\lambda)-1$, and if 
$D^\nu$ is another composition factor then 
$\epsilon_i(D^\nu)<\epsilon_i(D^\mu)$. 
\end{proposition}

In the following, we denote by $B$ the second crystal, and 
by $\mathcal{KP}$ the first crystal defined on the same set $B$. 

\section{Proof of the modular branching rule}

We assume the conditions $q\ne1$ and $v_i=q^{\gamma_i}$ as before. 

\begin{theorem}\label{main}
For $\lambda\in\mathcal{KP}$, we have that 
$\tilde e_iD^\lambda\ne0$ if and only if $\tilde e_i\lambda\ne0$ and if this 
holds then $\tilde e_iD^\lambda=D^{\tilde e_i\lambda}$. 
\end{theorem}
\begin{proof}
We first assume that the characteristic of $F$ is zero. 

As $\mathcal{KP}$ and $B=\{D^\lambda|\lambda\in\mathcal{KP}\}$ are 
isomorphic crystals by theorem \ref{Brundan-Kleshchev}, 
there exists a bijection 
$c:\mathcal{KP}\simeq \mathcal{KP}$ such that 
\begin{gather*}
\tilde e_iD^{c(\lambda)}=D^{c(\tilde e_i\lambda)},\quad
\tilde f_iD^{c(\lambda)}=D^{c(\tilde f_i\lambda)},\quad
wt(c(\lambda))=wt(D^{c(\lambda)})=wt(\lambda),\\[5pt]
\epsilon_i(D^{c(\lambda)})=\epsilon_i(\lambda),\quad
\varphi_i(D^{c(\lambda)})=\varphi_i(\lambda).
\end{gather*}
We prove by induction on $n$ that 
$c(\lambda)=\lambda$ for $\lambda\vdash n$. 
If $n=0$ there is nothing to prove. 
If $n=1$, $D^\lambda$ is the one dimensional module of the 
truncated polynomial 
ring $\H_1$ on which $L_1$ acts as $q^i\in\{v_1,\dots,v_m\}$ where 
$i$ is the color of the unique node of $\lambda$. 
Thus, $\tilde e_iD^\lambda=D^\emptyset=D^{c(\tilde e_i\lambda)}$ and 
$$
D^{c(\tilde e_ic^{-1}(\lambda))}=\tilde e_iD^{\lambda}=
D^{c(\tilde e_i\lambda)}\ne0.
$$
Then, $c(\tilde e_ic^{-1}(\lambda))=c(\tilde e_i\lambda)\ne0$, 
which implies $c(\lambda)=\lambda$. 

Assume that $n>1$ and that $c(\mu)=\mu$ for all $|\mu|<n$. 
Let $D^\mu=\tilde e_iD^\lambda\ne0$. Then, $c(\mu)=\mu$ implies
\begin{equation*}
\begin{split}
\epsilon_i(c^{-1}(\lambda))&=
\epsilon_i(D^\lambda)=\epsilon_i(D^\mu)+1=\epsilon_i(\mu)+1,\\[5pt]
\varphi_i(c^{-1}(\lambda))&=
\varphi_i(D^\lambda)=\varphi_i(D^\mu)-1=\varphi_i(\mu)-1.
\end{split}
\end{equation*}
By theorem \ref{Fock space theory} and lemma \ref{Kashiwara}, we have 
$$
f_iP^\mu=(P^{\tilde f_i\mu})^{\oplus(\epsilon_i(\mu)+1)}\bigoplus 
\left(\bigoplus_{\lambda'}(P^{\lambda'})^{\oplus a_{\mu\lambda'}^i}\right),
$$
where $a_{\mu\lambda'}^i$ are certain nonnegative integers, and 
$\lambda'$ satisfy $\lambda'\vdash n$ and 
$$
\epsilon_i(\lambda')\ge\epsilon_i(\mu)+2>\epsilon_i(c^{-1}(\lambda)).
$$
As $D^\lambda=\tilde f_iD^\mu=\Top(f_iD^\mu)$ and 
we have surjection
$f_iP^\mu \rightarrow f_iD^\mu$, $\lambda$ is either 
$\tilde f_i\mu$ or one of $\lambda'$. 
If $\lambda=\tilde f_i\mu$ then
$$
D^\mu=\tilde e_iD^{\tilde f_i\mu}
=D^{c(\tilde e_ic^{-1}(\tilde f_i\mu))}
=D^{\tilde e_ic^{-1}(\tilde f_i\mu)}. 
$$
Thus $\mu=\tilde e_ic^{-1}(\tilde f_i\mu)\ne0$ implies 
$\tilde f_i\mu=c^{-1}(\tilde f_i\mu)$ and $c(\lambda)=\lambda$ 
follows. Hence, we may assume 
$\epsilon_i(\lambda)>\epsilon_i(c^{-1}(\lambda))$. 
Next, we consider 
$$
e_iP^\lambda=(P^{\tilde e_i\lambda})^{\oplus(\varphi_i(\lambda)+1)}\bigoplus 
\left(\bigoplus_{\mu'}(P^{\mu'})^{\oplus b_{\lambda\mu'}^i}\right),
$$
where $b_{\lambda\mu'}^i$ are certain nonnegative integers, and 
$\mu'$ satisfy $\mu'\vdash n-1$ and 
$$
\varphi_i(\mu')\ge\varphi_i(\lambda)+2.
$$
Recall that $\H_n$ is a symmetric algebra. 
As $D^\mu=\tilde e_iD^\lambda=\Soc(e_iD^\lambda)$ and 
we have injection
$e_iD^\lambda \rightarrow e_iP^\lambda$, $\mu$ is either 
$\tilde e_i\lambda$ or one of $\mu'$. 
If $\mu=\tilde e_i\lambda$ then
$$
D^\lambda=\tilde f_iD^{\tilde e_i\lambda}
=D^{c(\tilde f_ic^{-1}(\tilde e_i\lambda))}
=D^{c(\tilde f_i\tilde e_i\lambda)}=D^{c(\lambda)}. 
$$
Thus $c(\lambda)=\lambda$ again follows. 
Hence, we may assume 
$\varphi_i(\mu)\ge\varphi_i(\lambda)+2$. 
As $\varphi_i(c^{-1}(\lambda))=\varphi_i(\mu)-1$, this implies 
$\varphi_i(c^{-1}(\lambda))>\varphi_i(\lambda)$. 

If both $\epsilon_i(\lambda)>\epsilon_i(c^{-1}(\lambda))$ and 
$\varphi_i(c^{-1}(\lambda))>\varphi_i(\lambda)$ hold, 
$$
\varphi_i(c^{-1}(\lambda))-\epsilon_i(c^{-1}(\lambda))
>\varphi_i(\lambda)-\epsilon_i(c^{-1}(\lambda))
>\varphi_i(\lambda)-\epsilon_i(\lambda).
$$
Thus $wt(c^{-1}(\lambda))(h_i)>wt(\lambda)(h_i)$, 
which contradicts to $wt(c^{-1}(\lambda))=wt(\lambda)$. 
We have proved the theorem when $F$ is of characteristic zero. 

Now we consider the positive characteristic case. Let $(K,R,F)$ 
be a modular system with parameters such that 
the characteristic of $K$ is zero, $\hat q\in R$ is a primitive $e^{th}$ 
root of unity, and $\hat q$ maps to $q\in F$. The image of $S_R^\lambda$ in 
$D_K^\lambda$ is denoted by $D_R^\lambda$. Since both $\hat q$ and 
$q$ have the multiplicative order $e$, we have 
$e_iD_K^\lambda=e_iD_R^\lambda\otimes_R K$. We also have surjection 
$e_iD_R^\lambda\otimes_R F\rightarrow e_iD_F^\lambda$, 
because $e_i$ is exact. Since we can read 
$\epsilon_i(D^\lambda)=
\max\{k\in\mathbb Z_{\ge0}|\tilde e_i^kD^\lambda\ne0\}$ from 
its restriction to the commutative subalgebra generated by 
$L_1,\dots,L_n$, that 
we have surjection $D_R^\lambda\rightarrow D_F^\lambda$ and 
injection $D_R^\lambda\rightarrow D_K^\lambda$ implies 
$\epsilon_i(D_K^\lambda)\ge\epsilon_i(D_F^\lambda)$. 
However, theorem \ref{Brundan-Kleshchev} guarantees that 
the sum of the left hand side and the right hand side in 
each weight space is the same. Hence, 
by the proof for the characteristic zero case, we have 
$\epsilon_i(D_F^\lambda)=\epsilon_i(D_K^\lambda)=\epsilon_i(\lambda)$ 
and, by proposition \ref{Grothendieck level}, 
$\tilde e_iD_F^\lambda=\Soc(e_iD_F^\lambda)$ is the 
composition factor $D_F^\mu$ of $e_iD_F^\lambda$ 
with the value $\epsilon_i(D_F^\mu)=\epsilon_i(\lambda)-1$. 
Observe that $e_iD_F^\lambda$ is self dual. Thus 
$\Top(e_iD_F^\lambda)=D_F^\mu$. 
Let $P_R^\mu$ be the lift of $P_F^\mu$. Then we have surjection 
$P_R^\mu\rightarrow e_iD_F^\lambda$. Consider the surjection 
$e_iD_R^\lambda\rightarrow e_iD_R^\lambda\otimes_R F
\rightarrow e_iD_F^\lambda$. 
Then the surjection $P_R^\mu\rightarrow e_iD_F^\lambda$ lifts to 
$P_R^\mu\rightarrow e_iD_R^\lambda$, which we denote by $f$. 
Recall that $P_R^\mu$ has 
Specht filtration $P_R^\mu=F_0\supset F_1\supset\cdots$ 
such that $F_0/F_1=S_R^\mu$. Let $f'$ be the 
composition of $f$ with the surjection 
$e_iD_R^\lambda\rightarrow \Top(e_iD_F^\lambda)$. 
As $f'$ factors through $P_F^\mu=P_R^\mu\otimes_R F$ and 
$F_1\otimes_R F$ is a proper submodule of $P_F^\mu=F_0\otimes_R F$ 
because $F_0\otimes_R F/F_1\otimes_R F=S_F^\mu$, 
that $P^\mu_F$ is the projective cover of 
$\Top(e_iD_F^\lambda)$ implies that 
$f'(F_1)=0\subset \Top(e_iD_F^\lambda)=f'(F_0)$. We have proved 
$f(F_0)\ne f(F_1)$. 
Let $K=\operatorname{Ker}\;f$. Then we have 
$$
0\rightarrow K \rightarrow F_0 \rightarrow f(F_0)\rightarrow 0.
$$
Since these are free $R$-modules, the exact sequence splits as 
$R$-modules. Thus, there exists a surjective $R$-linear map
$$
f(F_0)/f(F_1)\bigoplus K/K\cap F_1\rightarrow F_0/F_1=S_R^\mu.
$$
Suppose that $f(F_0)$ and $f(F_1)$ have the same rank as $R$-modules. 
Then $K/K\cap F_1\rightarrow F_0/F_1$ is surjective, since 
$f(F_0)/f(F_1)$ is a torsion $R$-module and $S_R^\mu$ is a free $R$-module. 
Thus, for any $x\in F_0$, we may write $x=y+z$ where 
$y\in F_1$ and $z\in K$, which implies $f(x)=f(y)\in f(F_1)$ 
and $f(F_0)=f(F_1)$, which is a contradiction. Therefore, 
we must have $f(F_0)\otimes_R K\ne f(F_1)\otimes_R K$. 
Consider the surjection 
$S_K^\mu\rightarrow f(F_0)\otimes_R K/f(F_1)\otimes_R K$. 
Since $\mu$ is Kleshchev and since 
$f(F_0)\otimes_R K/f(F_1)\otimes_R K\ne0$, 
the kernel of the map is contained in 
$\Rad S_K^\mu$. Thus we have 
$[f(F_0)\otimes_R K/f(F_1)\otimes_R K:D_K^\mu]\ne0$. 
As $f(F_0)\otimes_R K/f(F_1)\otimes_R K$ is a subquotient of 
$e_iD_R^\lambda\otimes_R K$, $D_K^\mu$ appears as a composition 
factor of $e_iD_K^\lambda$ with 
$\epsilon_i(D_K^\mu)\ge\epsilon_i(\lambda)-1$. As 
the maximum value in $e_iD_K^\lambda$ is 
$\epsilon_i(\lambda)-1$ and it is attained by 
$D_K^{\tilde e_i\lambda}$ by the proof in 
the characteristic zero case, we conclude that 
$\mu=\tilde e_i\lambda$ as desired. 
\end{proof}

\begin{remark}
As a corollary, $\operatorname{dim} D^\lambda$ is greater than 
or equal to the number of paths from $\emptyset$ to $\lambda$ 
in $\mathcal{KP}$. 
\end{remark}

\bibliographystyle{amsalpha}

\end{document}